\documentclass[12pt]{amsart}
\usepackage{amsmath,amscd,amssymb,amsfonts}
\newcommand{\h}{\hbox}
\newcommand{\q}{\quad}
\newcommand{\nin}{\par\noindent}
\newcommand{\bs}{\par\bigskip}
\newcommand{\ms}{\par\medskip}
\newcommand{\sk}{\par\smallskip}
\newcommand{\bsn}{\par\bigskip\noindent}
\newcommand{\msn}{\par\medskip\noindent}
\newcommand{\skn}{\par\smallskip\noindent}
\newcommand{\ssb}{\raise.15ex\h{${\scriptscriptstyle\bullet}$}}
\newcommand{\ssc}{\,\raise.15ex\hbox{${\scriptstyle\circ}$}\,}
\newcommand{\msum}{\h{$\sum$}}
\newcommand{\mopl}{\h{$\bigoplus$}}
\newcommand{\C}{{\mathbb C}}
\newcommand{\N}{{\mathbb N}}
\newcommand{\PP}{{\mathbb P}}
\newcommand{\Q}{{\mathbb Q}}
\newcommand{\R}{{\mathbf R}}
\newcommand{\Z}{{\mathbb Z}}
\newcommand{\B}{{\mathcal B}}
\newcommand{\D}{{\mathcal D}}
\newcommand{\F}{{\mathcal F}}
\newcommand{\G}{{\mathcal G}}
\newcommand{\Hc}{{\mathcal H}}
\newcommand{\OO}{{\mathcal O}}
\newcommand{\BB}{\widetilde{\mathcal B}}
\newcommand{\bt}{\widetilde{b}}
\newcommand{\m}{\widetilde{m}}
\newcommand{\M}{{}\,\widetilde{\!M}{}}
\newcommand{\Ht}{\widetilde{H}}
\newcommand{\dd}{{\partial}}
\newcommand{\ddd}{{\rm d}}
\newcommand{\ee}{{\mathbf e}}
\newcommand{\al}{\alpha}
\newcommand{\be}{\beta}
\newcommand{\la}{\lambda}
\newcommand{\om}{\omega}
\newcommand{\Om}{\Omega}
\newcommand{\De}{\Delta}
\newcommand{\Gr}{{\rm Gr}}
\newcommand{\DR}{{\rm DR}}
\newcommand{\End}{{\rm End}}
\newcommand{\bl}{\bigl}
\newcommand{\br}{\bigr}
\newcommand{\into}{\hookrightarrow}
\newcommand{\onto}{\mathop{\rlap{$\to$}\hskip2pt\hbox{$\to$}}}
\newcommand{\simto}{\,\,\rlap{\hskip1.3mm\raise1.4mm\hbox{$\sim$}}\hbox{$\longrightarrow$}\,\,}
\newcommand{\ges}{\geqslant}
\newcommand{\les}{\leqslant}
\begin{document}
\title[$D$-modules generated by powers of functions]
{$D$-modules generated by rational powers\\
of holomorphic functions}
\author[M. Saito]{Morihiko Saito}
\address{RIMS Kyoto University, Kyoto 606-8502 Japan}
\dedicatory{Dedicated to Professor Masaki Kashiwara}
\begin{abstract} We prove some sufficient conditions in order that a root of the Bernstein-Sato polynomial contributes to a difference between certain $D$-modules generated by rational powers of a holomorphic function; for instance, this holds in the case of isolated singularities with semisimple Milnor monodromies. We then construct an example where a root does not contribute to a difference. This also solves an old open problem about the relation between the Milnor monodromy and the exponential of the residue of the Gauss-Manin connection on the saturation of the Brieskorn lattice. This shows that the structure of Brieskorn lattices can be more complicated than one might imagine.
\end{abstract}
\maketitle
\centerline{\bf Introduction}
\bsn
Let $f$ be a holomorphic function on a complex manifold $X$. Set $Z:=f^{-1}(0)\subset X$. We assume $n:=\dim X\ges 2$ in this paper. For $\al\in\Q_{>0}$, consider the $\D_X$-modules
$$\D_Xf^{-\al}\subset\OO_X(*Z)f^{-\al},$$
where $\OO_X(*Z)=\OO_X\bl[\frac{1}{f}\br]$.
Let $b_f(s)$ be the Bernstein-Sato polynomial of $f$, see \cite{Be}, \cite{Bj}, \cite{Sat}, etc. It is the monic polynomial in $s$ satisfying the functional equation
$$b_f(s)f^s=Pf^{s+1}\q\h{with}\,\,\,\,P\in\D_X[s],$$
and having the minimal degree. Specializing this equation to $s=-\al$, we get
$$\D_Xf^{-\al}=\D_Xf^{-\al+1}\,\,\,\,\h{if}\,\,\,\,b_f(-\al)\ne 0.$$
\sk
One may then naturally ask the following
\msn
{\bf Question 1.} $\,\D_Xf^{-\al}\ne\D_Xf^{-\al+1}\,$ if $\,b_f(-\al)=0\,$?
\ms
Indeed, this has been asked to the author also by N.~Budur and U.~Walther recently.
In this paper we show that it can be solved to some extent (see Theorems~1, 2, and 3 below) by using the theory of Bernstein-Sato polynomials and $V$-filtrations as in \cite{Ka1}, \cite{Ka2}, \cite{Ka3}, \cite{Ma2}, \cite{Ma3}, \cite{mhp}, etc. (We may assume $\al\ne1$ in Question~1, since it is trivial for $\al=1$.)
It turns out, however, that there is an example where Question~1 has a negative answer, see Example~(4.2) below.
This also solves an old open problem in \cite{sup}, see Remark~(4.4) below.
\sk
Let $i_f:X\into Y:=X\times\C$ be the graph embedding by $f$ with $t$ the coordinate of $\C$. Set
$$M_f:=(i_f)_*^{\!\D}\OO_X,\q M_f^{\al}:=\Gr_V^{\al}M_f\,\,\,\,(\al\in\Q).$$
Here $(i_f)_*^{\!\D}\OO_X$ denotes the direct image of the structure sheaf $\OO_X$ {\it as $\D$-module,} and $V$ is the filtration of Kashiwara \cite{Ka3} and Malgrange \cite{Ma3} which is indexed by $\Q$ so that $\dd_tt-\al$ is locally nilpotent on $\Gr_V^{\al}M_f$, see also \cite{mhp}, \cite{rat}, etc.
Malgrange's formula \cite{Ma3} implies that
$$b_f(-\al)=0\iff\Gr^G_jM_f^{\be}\ne 0\q\,\,\,\bl(\al=\be+j,\,\,\be\in(0,1],\,\,j\in\N\br),
\leqno(1)$$
where $G$ is a filtration on $M_f^{\al}$ defined by using the inclusion $\D_X[s]f^s\into M_f$, see (1.3) below.
\sk
Set $N:=s+\be=-\dd_tt+\be$ on $M_f^{\be}$. In this paper we show the following
\msn
{\bf Theorem~1.\,} {\it There are canonical isomorphisms of regular holonomic $\D_X$-modules}
$$\D_Xf^{-\al}/\D_Xf^{-\al+1}=\Gr^G_j\bl(M_f^{\be}/NM_f^{\be}\br)\q\,\,\,\bl(\al=\be+j,\,\,\be\in(0,1],\,\,j\in\N\br).
\leqno(2)$$
\ms
Notice that there is a difference between $\Gr^G_jM_f^{\be}$ in (1) and $\Gr^G_j\bl(M_f^{\be}/NM_f^{\be}\br)$ in (2). This causes some problems.
Set
$$\M_f^{\be}:=\begin{cases}M_f^{\be}&\h{if}\,\,\,\be\in(0,1),\\ M_f^1/{\rm Ker}\,N&\h{if}\,\,\,\be=1.\end{cases}$$
This corresponds to the passage from the nearby cycles $\psi_f\C_X$ to the vanishing cycles $\varphi_f\C_X$, see the explanations after Theorems~2 and 3 below. We have the quotient filtration of $G$ defined by the surjection $M_f^{\be}\onto\M_f^{\be}$. It is denoted also by $G$.
\msn
{\bf Proposition~1.\,} {\it We have the following isomorphisms for $(\be,j)\ne(1,0)$\,$:$
$$\Gr^G_jM_f^{\be}\simto\Gr^G_j\M_f^{\be},\q\q\Gr^G_j\bl(M_f^{\be}/NM_f^{\be}\br)\simto\Gr^G_j\bl(\M_f^{\be}/N\M_f^{\be}\br),
\leqno(3)$$
and the following short exact sequences of regular holonomic $\D_X$-modules for $(\be,j)=(1,0)$\,$:$
$$\aligned0\to M'_Z\to\Gr^G_0M_f^1&\to\Gr^G_0\M_f^1\to 0,\\
0\to M_Z\to\Gr^G_0(M_f^1/NM_f^1)&\to\Gr^G_0\bl(\M_f^1/N\M_f^1\br)\to 0.\endaligned
\leqno(4)$$
Here $M'_Z$, $M_Z$ are the regular holonomic $\D_X$-modules such that $\DR_X(M'_Z)=\C_Z[n-1]$, $\DR_X(M_Z)=\h{\rm IC}_Z\C$ respectively with $\h{\rm IC}_Z\C$ the intersection complex of $Z\subset X$.}
\ms
This is shown in (1.4) below. (Note that $G_jM_f^{\be}=0$ if $j<0$.)
By Proposition~1, the assertions (1) and (2) hold with $M_f^{\be}$ replaced by $\M_f^{\be}$ if $\al\ne1$.
This is compatible with \cite{Ma2} in the isolated singularity case and with \cite{mic} in the general case.
\sk
From Theorem~1 and Proposition~1, we can deduce the following.
\msn
{\bf Theorem~2.\,} {\it Let $\be\in(0,1]$ with $\al=\be+j$ as above. We have a positive answer to Question~$1$ if $N$ is strictly compatible with the filtration $G$ on $\M_f^{\be}$, for instance, if $f$ has only an isolated singularity with semisimple Milnor monodromy $($that is, $N=0)$.}
\msn
{\bf Theorem~3.\,} {\it Set $\m_{\be}:=\min\{k\in\N\mid N^k\M_f^{\be}=0\}$. Assume $-\al$ is a root of $b_f(s)$ with the maximal multiplicity $\m_{\be}$, where $\al=\be+j$ as above. Then we have a positive answer to Question~$1$.}
\ms
The assumption on the strict compatibility of $N$ with $G$ in Theorem~2 is satisfied, for instance, if $N=0$ on $\M_f^{\be}$, or if $G$ is trivial on $\M_f^{\be}$ (that is, if there is no root $\al'$ of $\bt_f(s):=b_f(s)/(s+1)$ with $\al'-\al\in\Z\setminus\{0\}$, see \cite{mic}).
However, this strict compatibility condition is not always satisfied, and Question~1 may have a negative answer even in the isolated singularity case, see Example~(4.2) below.
\sk
We now explain the meaning of $\M_f^{\be}$. In this paper we denote by $\PP(X,\C)$ the abelian full subcategory of $D^b_c(X,\C)$ which is defined in \cite{BBD}, and is equivalent to the category of regular holonomic $\D_X$-modules via the de Rham functor $\DR_X$ by Kashiwara and Mebkhout as is well-known. By \cite{Ka3}, \cite{Ma3} we have the following isomorphisms for $\la=e^{-2\pi i\be}$, $\be\in(0,1]$\,:
$$\DR_X\bl(M_f^{\be}\br)=\psi_{f,\la}\C_X[n-1],\,\,\,\DR_X\bl(\M_f^{\be}\br)=\varphi_{f,\la}\C_X[n-1]\,\,\,\,\h{in}\,\,\,\PP(X,\C).
\leqno(5)$$
Here the de Rham complexes $\DR_X\bl(M_f^{\be}\br)$, $\DR_X\bl(\M_f^{\be}\br)$ are shifted by $n=\dim X$, and $\psi_{f,\la}$, $\varphi_{f,\la}$ respectively denote the $\la$-eigenpart of the nearby and vanishing cycle functors $\psi_f$, $\varphi_f$, see \cite{Dvc}, and also \cite{mhp}, etc. Moreover, $N=s+\be$ on the left-hand side of the equalities corresponds to $N:=(2\pi i)^{-1}\log T_u$ on the right-hand side, where $T_u$ is the unipotent part of the monodromy.
This implies that $M_f^{\be}$, $\M_f^{\be}$ in (1), (2), and Theorem~2 can be replaced respectively with $\psi_{f,\la}\C_X[n-1]$, $\varphi_{f,\la}\C_X[n-1]\in\PP(X,\C)$. (Note, however, that the vanishing of $N$ on $\varphi_{f,\la}\C_X[n-1]$ is not equivalent to that on the Milnor cohomology unless $f$ has only {\it isolated singularities}.)
\sk
Let $L_{\la^{-1}}$ be a local system of rank 1 on $U:=X\setminus f^{-1}(0)$ which is the pull-back of a local system of rank 1 on $\C^*$ with monodromy $\la^{-1}$ by $f$. This corresponds to $\D_U(f|_U)^{-\be}$ by the de Rham functor $\DR_U$. There are isomorphisms
$$\psi_{f,\la}\C_X[n-1]=\psi_{f,1}L_{\la^{-1}}[n-1]\q\h{in}\,\,\,\PP(X,\C),$$
and extensions of $L_{\la^{-1}}[n]$ in $\PP(X,\C)$ can be described by the Verdier-type extension theorem \cite{Ve}. The corresponding extensions of $\D_U(f|_U)^{-\be}$ as regular holonomic $\D_X$-modules are described by regular holonomic $\D_X$-submodules of $M_f^{\be}$ containing $NM_f^{\be}$ if we consider only extensions having no nontrivial submodules supported in $f^{-1}(0)$, see (1.9) below. The right-hand side of (2) is closely related to this. Note that $NM_f^{\be}\,(\subset M_f^{\be})$ corresponds to the intermediate direct image in the sense of \cite{BBD}.
\sk
If $f$ has an isolated singularity at $0\in Z$ and moreover ${\rm Sing}\,Z=\{0\}$, then the vanishing cycle sheaf $\varphi_{f,\la}\C_X[n-1]$ can be identified with the $\la$-eigenpart of the vanishing cohomology $H^{n-1}(F_{\!f,0},\C)_{\la}$ where $F_{\!f,0}$ denotes the Milnor fiber of $f$ around $0\in Z$.
In this case there are canonical isomorphisms
$$\bl(\M_f^{\be},G\br)=(i_0)_*^{\!\D}\bl(\Gr_V^{\be}\G_f,G\br)=(i_0)_*^{\!\D}\bl(H^{n-1}(F_{\!f,0},\C)_{\ee(-\be)},G\br)\q(\be\in(0,1]),
\leqno(6)$$
with $\ee(-\be):=e^{-2\pi i\be}$. Here $(i_0)_*^{\!\D}$ denotes the direct image as $\D$-modules by the closed embedding $i_0:\{0\}\into X$ (which is given by the tensor product with $\C[\dd_{x_1},\dots,\dd_{x_n}]$ over $\C$), and $\G_f$ is the Gauss-Manin system having the $V$-filtration of Kashiwara \cite{Ka3} and Malgrange \cite{Ma3} indexed by $\Q$, see \cite{Ph}, \cite{bri}, \cite{SS}, etc. The filtration $G$ on the middle term is defined by using the the saturation of the Brieskorn lattice $H''_f$ (see \cite{Br}) as in \cite{Ma1}, and this induces the filtration $G$ on the last term by using the last isomorphism of (6).
By Theorem~1 and Proposition~1 we get the following.
\msn
{\bf Corollary~1.\,} {\it If $f$ is a weighted homogeneous polynomial with an isolated singularity at $0$, then in the notation of $(6)$, we have the following isomorphisms for $\al\ne1$ and $p=[n-\al]$\,$:$
$$\D_Xf^{-\al}/\D_Xf^{-\al+1}=(i_0)_*^{\!\D}\,\Gr_F^pH^{n-1}(F_{\!f,0},\C)_{\ee(-\al)},
\leqno(7)$$
and the following short exact sequence of regular holonomic $\D_X$-modules for $\al=1$\,$:$}
$$0\to M_Z\to\D_Xf^{-1}\!/\OO_X\to(i_0)_*^{\!\D}\,F^{n-1}H^{n-1}(F_{\!f,0},\C)_1\to 0.
\leqno(8)$$
\ms
Indeed, the filtration $G$ coincides with the Hodge filtration $F$ of the mixed Hodge structure on the Milnor cohomology (\cite{St1},\cite{St2}) up to a shift of filtration in this case, and there are well-known canonical isomorphisms
$$F^pH^{n-1}(F_{\!f,0},\C)_{\ee(-\al)}=\Gr_V^{\al}\Om_f\q(\al\in\Q,\,\,p=[n-\al]),
\leqno(9)$$
where $\Om_f:=\Om_{X,0}^n/\ddd f\wedge\Om_{X,0}^{n-1}$, see \cite{St1}, \cite{SS}, etc.
Corollary~1 seems to be closely related to a question of T.~Bitoun who has calculated with T.~Schedler the length of the regular holonomic $\D_X$-module $\D_Xf^{-1}\!/\OO_X$ in the case where $f$ defines a nonsingular curve $C$ in $\PP^2$. Note that $F^2H^2(F_{\!f,0},\C)_1=F^1H^1(C,\C)$ in this case. (They have recently generalized their result to the weighted homogenous isolated singularity case \cite{BiSc}.)
\sk
N.~Budur \cite{Bu} informed us that Question~1 has a positive answer if there is no root $-\al'$ of $b_f(s)$ with $\al'-\al\in\Z_{>0}$ (by using a different method, see also \cite{Oa}). This can be deduced also from Theorem~1 in this paper (since $\Gr^G_j\bl(M_f^{\be}/NM_f^{\be}\br)\ne 0$ if $\Gr^G_jM_f^{\be}\ne 0$ and $G_jM_f^{\be}=M_f^{\be}$).
\sk
I thank N. Budur and U. Walther for their interest in this problem, and T.~Bitoun for a good question related to Corollary~1. I also thank the referee for useful comments.
This work was partially supported by JSPS Kakenhi 15K04816.
\sk
In Section~1 we review some basics around the theory of Bernstein-Sato polynomials.
In Section~2 we prove Proposition~(2.3) which is essential for the proofs of the main theorems.
In Section~3 we prove the main theorems.
In Section~4 we describe an example where Question~1 is negatively answered.
\bs\bs
\vbox{\centerline{\bf \S1. Preliminaries}
\bsn
In this section we review some basics around the theory of Bernstein-Sato polynomials.}
\msn
{\bf 1.1.~Direct image by the graph embedding.\,} Let $i_f:X\into Y:=X\times\C$ be the graph embedding by a holomorphic function $f$ on a complex manifold $X$ with $n:=\dim X\ges 2$. Let $t$ be the coordinate of $\C$. We denote by $(i_f)_*^{\!\D}\OO_X$ the direct image of the structure sheaf $\OO_X$ as $\D$-module. It is a free $\OO_X[\dd_t]$-module of rank 1, and has the formal generator $\delta(f-t)$ satisfying the relations
$$\dd_{x_i}\,\delta(f-t)=-(\dd_{x_i}f)\,\dd_t\,\delta(f-t),\q t\,\delta(f-t)=f\,\delta(f-t),\,\,\h{etc.,}
\leqno(1.1.1)$$
with $x_1,\dots,x_n$ local coordinates of $X$.
This expresses the twist of the action of the vector fields $\dd_{x_i}$, see also \cite[Section 1.8]{rat}.
\sk
More precisely, for local coordinates $x_1,\dots,x_n$ of $X$, set
$$x'_i:=x_i,\q t':=t-f.$$
Then, for a $\D_X$-module $M$, there is a canonical isomorphism of $\D_Y$-modules
$$(i_f)^{\!\D}_*M=M[\dd_{t'}],
\leqno(1.1.2)$$
where the action of $\dd_{x'_i}$ on $M[\dd_{t'}]$ is the natural one, and
$$\dd_{x_i}=\dd_{x'_i}-(\dd_{x_i}f)\,\dd_{t'},\q\dd_t=\dd_{t'}.
\leqno(1.1.3)$$
Here the first equality, which corresponds to the first equality of (1.1.1), means the twist of the action of the vector fields $\dd_{x_i}$.
\sk
It is well-known (see \cite{Ma2}) that $\delta(t-f)$ and $f^s$ satisfy the same relations in $\D_Y$, and we have a natural inclusion
$$\D_X[s]f^s\into M_f:=(i_f)_*^{\!\D}\OO_X,
\leqno(1.1.4)$$
such that $f^s$ is naturally identified with $\delta(t-f)$, and $s$ with $-\dd_tt$.
(Here the sheaf-theoretic direct image by $i_f$ is omitted to simplify the notation.)
\sk
The roots of $b_f(s)$ are negative by \cite{Ka1}, and greater than $-n$ by \cite{mic}. (Recall that we assume $n\ges 2$ in this paper.) These respectively imply the inclusions
$$V^{>0}M_f\supset\D_X[s]f^s\supset V^{n-1}M_f,
\leqno(1.1.5)$$
where $V$ is the filtration of Kashiwara and Malgrange as in (1.2) below. The last inclusion can be shown by using the increasing filtration $G$ on $M_f\subset M_f[t^{-1}]$ in (1.3) below. Indeed, we can deduce the following from the above assertion in {\it loc.\,cit.}:
$$\Gr^G_jV^{\al}M_f=0\q\h{if}\q\al+j\ges n.
\leqno(1.1.6)$$
Here we may assume $j=0$ replacing $M_f$ with $M_f[t^{-1}]$ and using (1.3.4) below. Then (1.1.6) follows from {\it loc.\,cit.} This induces the desired inclusion by setting $\al=n-1$ in (1.1.6), since it means that $G_jV^{n-1}M_f=G_{j-1}V^{n-1}M_f$ for $j\ges 1$. (This is known for Brieskorn lattices in the hypersurface isolated singularity case.)
\msn
{\bf 1.2.~$V$-filtration.\,}
Let $M$ be a regular holonomic $\D_Y$-module, where $Y=X\times\C$ as in (1.1). We have the unique exhaustive filtration $V$ of Kashiwara \cite{Ka3} and Malgrange \cite{Ma3} which is indexed by $\Q$ and satisfies the following conditions:
\msn
\vbox{\nin
(a) $\,\,t\,(V^{\al}M)\subset V^{\al+1}M$ (with equality for $\al\gg 0$),\,\,\,$\dd_t\,(V^{\al}M)\subset V^{\al-1}M$\,\,\,($\al\in\Q$),
\msn
(b) $\,\,\dd_tt-\al$ is locally nilpotent on $\Gr_V^{\al}M$ ($\al\in\Q$),
\msn
(c) $\,\,V^{\al}M$ ($\forall\,\al\in\Q$) is locally finitely generated over $V^0\D_Y$.}
\msn
Here $V^0\D_Y\subset\D_Y$ is the subring generated by $\OO_Y$ and $\D_X[s]$ with $s:=-\dd_tt$.
(We can replace $V^0\D_Y$ with $\D_X[s]$ if $M$ is the direct image of a regular holomorphic $\D_X$-module by the graph embedding $i_f$.)
The existence of the $V$-filtration follows from that of $b$-functions in \cite{Ka2}, where the regularity is actually unnecessary, see also \cite[Proposition~1.9]{rat}.
Note that the increasing filtration $V_{\ssb}$ is used in \cite{mhp}, and we have $V_{\al}=V^{-\al}$.
\sk
The above conditions imply the isomorphisms
$$t:\Gr_V^{\al}M\simto\Gr_V^{\al+1}M,\q\dd_t:\Gr_V^{\al+1}M\simto\Gr_V^{\al}M\q(\al\ne 0).
\leqno(1.2.1)$$
Note that the first isomorphism holds for any $\al$ if and only if the action of $t$ on $M$ is bijective. (This follows from the Verdier-type extension theorem \cite{Ve} for regular holonomic $\D$-modules.)
\sk
It is known that the functor associating $V^{\al}M$ (or $\Gr_V^{\al}M$) to $M$ is exact, and
$$\h{$V^{\al}M=0$ for $\al>0$ if $M$ is supported in $\{t=0\}$,}
\leqno(1.2.2)$$
(see, for instance, \cite[Lemma 3.1.3]{mhp}).
These imply the following:
\sk
Let $M\to M'$ be a morphism of regular holonomic $\D_Y$-modules as above inducing an isomorphism over $\{t\ne 0\}$. Then it induces the isomorphisms
$$V^{\al}M\simto V^{\al}M'\q(\forall\,\al>0).
\leqno(1.2.3)$$
\msn
{\bf 1.3.~Filtration $G$.\,}
We have the increasing filtration $G$ on $M_f[t^{-1}]$ defined by
$$G_j(M_f[t^{-1}])=t^{-j}\D_X[s]f^s\subset M_f[t^{-1}]\q\,\,\,(j\in\Z).
\leqno(1.3.1)$$
This induces the increasing filtration $G$ on $M_f\subset M_f[t^{-1}]$ and on $M_f^{\be}=\Gr_V^{\be}M_f$ so that
$$G_jM_f^{\be}=t^{-j}\Gr_V^{\be+j}(\D_X[s]f^s)\,\subset\,M_f^{\be}\q\,\,\,\bl(\be\in(0,1],\,\,j\ges 0\br),
\leqno(1.3.2)$$
with $G_jM_f^{\be}=0$ ($j<0$) by the first inclusion of (1.1.5) (or \cite{Ka1}).
Here we have the following isomorphisms (using (1.2.3)):
$$M_f^{\al}\,(=\Gr_V^{\al}M_f)\simto Gr_V^{\al}(M_f[t^{-1}])\q\q\bl(\al>0),
\leqno(1.3.3)$$
$$t^k:(M_f[t^{-1}];V,G)\simto(M_f[t^{-1}];V[k],G[k])\q\q\bl(k\in\Z\br),
\leqno(1.3.4)$$
with $(V[k])^{\al}:=V^{\al+k}$, $(G[k])_j:=G_{j-k}$ (in a compatible way with $G_j=G^{-j}$). Note that $G$ is an increasing filtration since
$$t(\D_X[s]f^s)=\D_X[s]f^{s+1}\subset\D_X[s]f^s.
\leqno(1.3.5)$$
The formula (1) in the introduction then follows from the definition of Bernstein-Sato polynomial. Indeed, the latter says that $b_f(s)$ is the minimal polynomial of the action of $s$ on
$$\D_X[s]f^s/t(\D_X[s]f^s).
\leqno(1.3.6)$$
\sk
The multiplicity $m_{b,\al}$ of a root $-\al$ of $b_f(s)$ is given by
$$m_{b,\al}=\min\bl\{k\in\N\mid N^k\Gr^G_jM_f^{\be}=0\br\}\q(\al=\be+j,\,\,\be\in(0.1],\,\,j\in\N),
\leqno(1.3.7)$$
where $N=s+\be$ on $M_f^{\be}$ as in Theorem~1. Note that the $G_jM_f^{\be}$ are stable by the action of $s$ and $N$.
We may replace $\Gr^G_jM_f^{\be}$ with $\Gr^G_j\M_f^{\be}$ in (1.3.7) by Proposition~1 if $\al\ne 1$.
This is closely related with \cite{Ma2}, \cite{mic}.
\sk
In the isolated singularity case, (1.3.7) is also related with the Steenbrink spectrum \cite{St2} by \cite{Ma2}, \cite{SS}, \cite{Va}, see (1.7) below.
\msn
{\bf 1.4.~Microlocalization.\,} In the notation of (1.1), set
$$\B_f:=(i_f)^{\!\D}_*\OO_X\buildrel{\iota}\over\into\BB_f:=\B_f[\dd_t^{-1}],$$
where the last term is the algebraic microlocalization, see \cite{mic}.
We have the microlocal $V$-filtration on $\BB_f$ so that, by setting
$$M_f:=\Gr_V^{\al}\B_f,\q\M_f:=\Gr_V^{\al}\BB_f\q(\al\in\Q),$$
there are isomorphisms as follows (see {\it loc.~cit.}):
$$\aligned\Gr_V^{\al}\iota:M_f^{\al}&\simto\M_f^{\al}\q(\al\in\Q\setminus\Z_{>0}),\\
\dd_t:\M_f^{\al}&\simto\M_f^{\al-1}\q(\al\in\Q).\endaligned
\leqno(1.4.1)$$
These imply that $\Gr_V^1\iota$ induces the isomorphism
$$M_f^1/{\rm Ker}\,N\simto\M_f^1,
\leqno(1.4.2)$$
which justifies the definition of $\M_f^1$ in the introduction.
Indeed, (1.4.2) follows from the short exact sequence
$$0\to{\rm Ker}\,N\into M_f^1\buildrel{\dd_t}\over\to M_f^0\to0,
\leqno(1.4.3)$$
where we use the injectivity of
$$t:M_f^0\to M_f^1.$$
This injectivity and the surjectivity of
$$\dd_t:M_f^1\to M_f^0$$
follow from the property of $\C_X[n]\in\PP(X,\C)$ that it has no nontrivial sub nor quotient object supported in $Z$, see also \cite[Lemma 5.1.4]{mhp}, etc. (Here $\PP(X,\C)$ is as in the introduction.) Note that (1.4.3) corresponds by the de Rham functor $\DR_X$ to the following short exact sequence in $\PP(X,\C)$\,:
$$0\to\C_Z[n-1]\to\psi_{f,1}\C_X[n-1]\to\varphi_{f,1}\C_X[n-1]\to 0.
\leqno(1.4.4)$$
Here $\C_Z[n-1]\in\PP(X,\C)$, since $Z$ is a hypersurface in a smooth manifold $X$.
\sk
We have moreover
$${\rm Ker}\,N\subset G_0M_f^1,
\leqno(1.4.5)$$
or equivalently
$$\C_Z[n-1]\subset G_0(\psi_{f,1}\C_X[n-1]),
\leqno(1.4.6)$$
where $G$ denotes also the corresponding filtration on $\psi_{f,1}\C_X[n-1]$.
\sk
We will show (1.4.6) by using the assertion that for $\F\in\PP(X,\C)$, we have
$${\rm Hom}(\C_Z[n-1],\F)=0\q\h{if}\q{\rm supp}\,\F\subset\Sigma:={\rm Sing}\,Z_{\rm red}.
\leqno(1.4.7)$$
Indeed, (1.4.7) follows from the adjunction isomorphism for the inclusion $i_{\Sigma}:\Sigma\into X$\,:
$$\R{\rm Hom}(\C_Z[\dim Z],\F)=\R{\rm Hom}(\C_{\Sigma}[\dim\Sigma],\F'[\dim\Sigma-\dim Z]),$$
where $\F'\in\PP(\Sigma,\C)$ with $\F=(i_{\Sigma})_*\F'$. We also use here the vanishing of negative extension groups in $\PP(X,\C)$ (see \cite{BBD}). Note that the sheaf complex $\C_{\Sigma}[\dim\Sigma]$ satisfies the condition
$$^p\Hc^j(\C_{\Sigma}[\dim\Sigma])=0\q(j>0),$$
by the definition of $\PP(\Sigma,\C)$, where $^p\Hc^j:D^b_c(X,\C)\to\PP(X,\C)$ is the cohomological functor defined in \cite{BBD}.
\sk
We apply (1.4.7) to the morphism
$$\C_Z[n-1]\to\psi_{f,1}\C_X[n-1]/G_0(\psi_{f,1}\C_X[n-1]),$$
where the target is supported in $\Sigma$. Indeed, setting $Z^{\circ}:=Z\setminus\Sigma$, we have
$$\C_{Z^{\circ}}[n-1]=\psi_{f,1}\C_X[n-1]|_{Z^{\circ}}=G_0(\psi_{f,1}\C_X[n-1])|_{Z^{\circ}}.$$
So (1.4.6) follows.
\sk
We can now prove Proposition~1 by using (1.4.5). Here we may assume $\be=1$, and we have the isomorphism
$$M'_Z={\rm Ker}\,N\subset M_f^1,$$
corresponding to a well-known isomorphism
$$\C_Z[n-1]={\rm Ker}\,N\subset\psi_{f,1}\C_X[n-1],$$
which follows from (1.4.4).
\sk
Moreover the kernel of the morphism
$$M_f^1/NM_f^1\to\M_f^1/N\M_f^1=M_f^1/({\rm Ker}\,N+{\rm Im}\,N)$$
is given by the image of ${\rm Ker}\,N\subset M_f^1$ in ${\rm Coker}\,N$, and the latter coincides with $M_Z$.
The last assertion follows from the $N$-primitive decomposition of $\Gr^W_{\ssb}\psi_{f,1}\C_X[n-1]$ together with the isomorphism
$$\Gr^W_{n-1}(\C_X[n-1])={\rm IC}_Z\C.
,$$
which also follows from (1.4.7) by using the semisimplicity of pure Hodge modules.
\sk
To conclude the proof of Proposition~1, note that the filtrations $G$ on $\M_f^1$, $M_f^1/NM_f^1$, and $\M_f^1/N\M_f^1$ are all quotient filtrations  of $G$ on $M_f^1$. By (1.4.5) we have the following isomorphisms for any $j\ges 0$\,:
$$\aligned M_f^1/G_jM_f^1&\simto\M_f^1/G_j\M_f^1,\\ (M_f^1/NM_f^1)\big/G_j(M_f^1/NM_f^1)&\simto\bl(\M_f^1/N\M_f^1\br)\big/G_j\bl(\M_f^1/N\M_f^1\br).\endaligned$$
By combining these with the above assertions about $M_Z$ and $M'_Z$, Proposition~1 then follows (since $G_jM_f^1=0$ for $j<0$).
\msn
{\bf 1.5.~Isolated singularity case.\,} Assume $Z$ has an isolated singularity at $0\in Z$. We have the Brieskorn lattice (see \cite{Br})
$$H''_f:=\Om_{X,0}^n/\ddd f\wedge\ddd\Om_{X,0}^{n-2},$$
which is contained in the Gauss-Manin system $\G_f$ (see \cite{Ph}). The latter is the localization of $H''_f$ by the action of $\dd_t^{-1}$ which is well-defined on $H''_f$. Indeed, we have
$$\dd_t^{-1}[\om]=[\ddd f\wedge\eta]\q\h{for}\,\,\,\,\om\in\Om_{X,0}^n,\,\,\eta\in\Om_{X,0}^{n-1}\,\,\,\,\h{with}\,\,\,\,\ddd\eta=\om.
\leqno(1.5.1)$$
It is known that $H''_f$ and $\G_f$ are respectively free modules of rank $\mu$ over $\C\{\!\{\dd_t^{-1}\}\!\}$ and $\C\{\!\{\dd_t^{-1}\}\!\}[\dd_t]$, see \cite{bri}, \cite{SS}, etc.
\sk
By Malgrange \cite{Ma2}, we have
$$\bt_f(s):=b_f(s)/(s+1)={\rm Min.\,Pol}(s,\Ht''_f/t\Ht''_f),
\leqno(1.5.2)$$
where the left-and side is the microlocal (or reduced) $b$-function (see \cite{mic}, \cite{Sat}), and the right-hand side is the minimal polynomial of the action of
$$s:=-\dd_tt\q\h{on}\q\Ht''_f/t\Ht''_f.$$
Here $\Ht''_f$ is the saturation of $H''_f$, that is,
$$H''_f\subset\Ht''_f:=\msum_{k\ges 0}\,s^kH''_f\subset\G_f.
\leqno(1.5.3)$$
This is closely related with the microlocalization in (1.4), and is generalized in some sense to the non-isolated singularity case in \cite{mic}.
\msn
{\bf 1.6.~Thom-Sebastiani type theorems.\,} Let $f,g$ be holomorphic functions with isolated singularities on complex manifolds $X,Y$. We denote by $f+g$ the function on $X\times Y$ defined by the sum of the pull-backs of $f$ and $g$. We have the Thom-Sebastiani type theorem for Brieskorn lattices (see \cite{Ma1} \cite{SS})
$$H''_f\otimes_{\C\{\!\{\dd_t^{-1}\}\!\}}H''_g\simto H''_{f+g}.
\leqno(1.6.1)$$
Indeed, there is a well-defined morphism with source and target free modules of the same rank over $\C\{\!\{\dd_t^{-1}\}\!\}$, and it induces an isomorphism by taking mod $\dd_t^{-1}\C\{\!\{\dd_t^{-1}\}\!\}$.
Note that the action of $t$ on the right-hand side (which is induced by the multiplication by $f+g$) corresponds to the action of $t\otimes id+id\otimes t$ on the left-hand side.
\sk
We then get the Thom-Sebastiani type theorem for Gauss-Manin systems
$$\G_f\otimes_K\G_g\simto\G_{f+g},
\leqno(1.6.2)$$
where $K:=\C\{\!\{\dd_t^{-1}\}\!\}[\dd_t]$.
Note that the action of $\dd_tt$ on the right-hand side corresponds to the action of $\dd_tt\otimes id+id\otimes\dd_tt$ on the left-hand side.
\msn
{\bf 1.7.~Spectral numbers.\,} Let
$$0<\al_{f,1}\les\cdots\les\al_{f,\mu}<n$$
be the spectral numbers of a holomorphic function $(f,0)$ with an isolated singularity (see \cite{St2}, \cite{Ar}, etc.), where $\mu$ is the Milnor number. These are defined by using the mixed Hodge structure on the vanishing cohomology, that is
$$\#\{i\mid\al_{f,k}=\al\}=\dim\Gr_F^pH_{\la}\q\q\bl(\la=\exp(-2\pi i\al),\,\,p=[n-\al]\br),$$
where $F$ is the Hodge filtration on the $\la$-eigenpart of the vanishing cohomology
$$H_{\la}:=H^{n-1}(F_{\!f,0},\C)_{\la},$$
in the notation of the introduction. We have moreover the following isomorphisms
$$\Gr_F^pH_{\la}=\Gr_V^{\al}\Om^n_f,
\leqno(1.7.1)$$
for $\la=\exp(-2\pi i\al)$, $p=[n-\al]$, see \cite{SS}, \cite{Va}. Here
$$\Om^n_f:=\Om_{X,0}^n/\ddd f\wedge\Om_{X,0}^{n-1}\cong H''_f/\dd_t^{-1}H''_f,$$
which is isomorphic to the Jacobian ring $\OO_{X,0}/(\dd f)$, and the $V$-filtration on $\Om^n_f$ is defined by using the last isomorphism.
\sk
We have the Thom-Sebastiani type theorem \cite{ST} for spectrum
$$\{\al_{f+g,k}\}_k=\{\al_{f,i}+\al_{g,j}\}_{i,j},
\leqno(1.7.2)$$
where $f,g$ are as in (1.6), see \cite{SS}, \cite{Va}.
\msn
{\bf 1.8.~Comparison with the spectral numbers.\,} We denote the eigenvalues of the action of $\dd_tt=-s$ on $\Ht''_f/t\Ht''_f$ (counted with multiplicities) by
$$0<\be_{f,1}\les\cdots\les\be_{f,\mu}<n.$$
These give the roots of $\bt_f(s)$ up to a sign by (1.5.2) by forgetting the multiplicities. (However, it is quite non-trivial to determine the latter.)
\sk
It follows from (1.7.1) together with the first inclusion of (1.5.3) that there is a permutation $\sigma$ of $\{1,\dots,\mu\}$ such that
$$\delta_{f,k}:=\al_{f,k}-\be_{f,\sigma(k)}\in\Z_{\ges0}.
\leqno(1.8.1)$$
This gives a certain estimate of the roots of $\bt_f(s)$.
However, $\sigma$ and the $\delta_{f,k}$ are not unique, and it is rather complicated to determine explicitly these in general (except for the case of weighted homogenous polynomials where $H''_f=\Ht''_f$ and hence $\al_{f,k}=\be_{f,k}$).
\sk
We have a better estimation by using a basis $(v_{r,i})_i$ of the $N$-primitive part
$$P\Gr^W_r\Gr_F^pH_{\la}:={\rm Ker}\,N^{r+1}\subset\Gr^W_r\Gr_F^pH_{\la}$$
where $W$ is the monodromy filtration with center 0. (So there is a shift of filtration with the one in \cite{St2} by $n-1$ if $\la\ne1$, and $n$ if $\la=1$.) By the primitive decomposition
$$\Gr^W_w\Gr_F^pH_{\la}=\mopl_{r\ges\max(0,-w)}\,N^rP\Gr^W_{w+2r}\Gr_F^{p+r}H_{\la},
\leqno(1.8.2)$$
we get a basis of $\Gr^W_w\Gr_F^pH_{\la}$ compatible with the above decomposition.
\sk
This implies that each spectral number $\al_{f,k}$ is associated with two additional integers $w_{f,k}$ and $r_{f,k}$, which are respectively given by $w$ and $r$ in (1.8.2).
There is no direct relation between $w_{f,k}$ and $r_{f,k}$.
Indeed, $w_{f,k}$ is used in the definition of spectral pairs (see \cite{St2}), but $r_{f,k}$ has not been considered before.
\sk
By using the primitive decomposition (1.8.2) together with (1.5.3) and (1.7.1), we get a better estimate
$$\delta'_{f,k}:=\al_{f,k}-r_{f,k}-\be_{f,\sigma(k)}\in\Z_{\ges0},
\leqno(1.8.3)$$
although it is more difficult to determine $\sigma$ in general. (Note that (1.8.3) is closely related to the fact that the Brieskorn lattice is a $(B)$-lattice in the sense of \cite{bri}.)
\msn
{\bf 1.9.~Verdier-type extension theorem.\,}
Let $Y=X\times\C$ as in (1.1). Set $Y_0:=\{t=0\}\subset Y$. Let $M$ be a regular holonomic $\D_Y$-modules such that the action of $t$ is bijective. Then there is a bijection between
$$\bl\{\h{Regular holonomic $\D_Y$-submodules $M'\subset M$ with ${\rm supp}\,M/M'\subset Y_0$}\br\}$$
and
$$\bl\{\h{Regular holonomic $\D_X$-submodules $M''\subset\Gr_V^0M$ containing $N\Gr_V^0M$}\}.$$
Here ``regular holonomic" can be replaced with ``coherent", since these are equivalent for $\D$-submodules of regular holonomic $\D$-modules. The correspondence is given by
$$M'\mapsto M'':=\Gr_V^0M'\subset\Gr_V^0M.
\leqno(1.9.1)$$
Indeed, let $M^{\rm min}\subset M$ be the minimal $\D_Y$-submodule with ${\rm supp}\,M/M^{\rm min}\subset Y_0$. This corresponds to the intermediate direct image in \cite{BBD}, and we have
$$\Gr_V^0M^{\rm min}={\rm Im}\,N\subset \Gr_V^0M,$$
see \cite[Proposition 3.1.8]{mhp}.
Moreover the above $M'$ (satisfying the above condition) contain $M^{\rm min}$, and are determined by the $\D_Y$-submodules
$$M'/M^{\rm min}\subset M/M^{\rm min}.$$
Indeed, it is enough to take the inverse image under the projection
$$M\onto M/M^{\rm min}.$$
Since $M'/M^{\rm min}$ is supported in $Y_0$, it is determined by $\Gr_V^0$, see \cite[Lemma 3.1.3]{mhp}. Moreover the functor $M'\mapsto\Gr_V^0M'$ is exact. So the assertion follows.
(The corresponding assertion in $\PP(X,\C)$ follows from \cite{Ve}. See also \cite{SZ} for the 1-dimensional case.)
\bs\bs
\vbox{\centerline{\bf \S2. Proof of the key proposition}
\bsn
In this section we prove Proposition~(2.3) which is essential for the proofs of the main theorems.}
\msn
{\bf 2.1.~Proposition.\,} {\it For $\al\in\Q$, there are canonical isomorphisms and an inclusion of regular holonomic $\D_Y$-modules
$$(i_f)_*^{\!\D}(\D_Xf^{-\al})=\D_Yf^{s-\al}=\D_X[\dd_t]f^{s-\al}\,\subset\,\OO_X(*Z)[s]f^{s-\al},
\leqno(2.1.1)$$
where the action of $t$ is bijective on the last $\D_Y$-module.}
\msn
{\it Proof.\,} Since $\OO_X(*Z)f^{-\al}$ is the pull-back of $\OO_{\C}\bl[\frac{1}{t}\br]t^{-\al}$ as $\OO$-module with an integrable connection, it is a regular holonomic $\D_X$-module.
This implies that the $\D_X$-submodule
$$\D_Xf^{-\al}\subset\OO_X(*Z)f^{-\al}$$
is also regular holonomic.
We can verify that the action of $t$ on $\OO_X(*Z)[s]f^{s-\al}$ is bijective (although we can also use the Riemann-Hilbert correspondence for this).
The last $\D_Y$-module in (2.1.1) is then the localization of $\D_Yf^{s-\al}$ along $Z$, and is regular holonomic, if we admit the isomorphisms in (2.1.1). So it is enough to show the first isomorphism of (2.1.1), since the second isomorphism is easy.
\sk
We verify that
$$f^{-\al}\otimes\delta(f-t)\,\in\,(i_f)_*^{\!\D}(\D_Xf^{-\al})\q\h{and}\q f^{s-\al}\,\in\,\OO_X(*Z)[s]f^{s-\al}$$
satisfy the same relations in $\D_Y(*Z)\,(:=\D_Y[1/f])$. Here $\delta(t-f)$ is the formal generator associated with the direct image by $i_f$ as in (1.1), and is used to express the twist of the action of vector fields as in (1.1.3). So the assertion follows.
\msn
{\bf 2.2.~Proposition.\,} {\it For $\al>0$, we have the isomorphisms}
$$\aligned\D_Yf^{s-\al}=\D_X[\dd_t]f^{s-\al}&=\msum_{i\ges 0}\,\dd_t^i(\D_X[s]f^{s-\al})\\&=\msum_{i\ges 0}\,\dd_t^i\,V^0(\D_X[s]f^{s-\al}).\endaligned
\leqno(2.2.1)$$
\msn
{\it Proof.\,} Only the last equality is nontrivial. It is reduced to
$$\D_X[s]f^{s-\al}\subset\msum_{i\ges 0}\,\dd_t^i\,V^0(\D_X[s]f^{s-\al}).
\leqno(2.2.2)$$
By decreasing induction on $-\al<\be<0$, we prove
$$V^{\be}\D_X[s]f^{s-\al}\subset\msum_{i\ges 0}\,\dd_t^i\,V^0(\D_X[s]f^{s-\al}).
\leqno(2.2.3)$$
Using the isomorphisms in (1.2.1), this can be reduced to the bijectivity of
$$\dd_tt\in\End\bl(\Gr_V^{\be}\D_X[s]f^{s-\al}\br)\q(\be<0).$$
So the assertion follows. This finishes the proof of Proposition~(2.2).
\msn
{\bf 2.3.~Proposition.\,} {\it For $\al>0$, there is an equality in $\Gr_V^0(\OO_X(*Z)[s]f^{s-\al})$}
$$\Gr_V^0(\D_Yf^{s-\al})=\Gr_V^0(\D_X[s]f^{s-\al})+\dd_t\,\Gr_V^1(\D_Yf^{s-\al}).
\leqno(2.3.1)$$
\msn
{\it Proof.\,} Since the inclusion $\supset$ is clear, it is enough to show the inclusion $\subset$ for the proof of (2.3.1). Set
$$M:=\OO_X(*Z)[s]f^{s-\al},$$
By the last inclusion in Proposition~(2.1) together with the isomorphism (1.2.3), we get the inclusion
$$\Gr_V^0(\D_Yf^{s-\al})\subset\Gr_V^0M,
\leqno(2.3.2)$$
and the isomorphisms
$$\Gr_V^1(\D_Yf^{s-\al})=\Gr_V^1M=t\,(\Gr_V^0M).
\leqno(2.3.3)$$
\sk
By the Verdier-type extension theorem for regular holonomic $\D$-modules (see (1.9)), there is a unique regular holonomic $\D_Y$-submodule $M'\subset M$ such that $M/M'$ is supported in $\{t=0\}$ and moreover
$$\Gr_V^0M'=\Gr_V^0(\D_X[s]f^{s-\al})+\dd_t\,\Gr_V^1(\D_Yf^{s-\al})\q\h{in}\,\,\,\Gr_V^0M,
\leqno(2.3.4)$$
where (2.3.3) is used.
The assertion is then reduced to the inclusion of $\D_Y$-submodules
$$\D_Yf^{s-\al}\subset M'\q\h{in}\,\,\,M.
\leqno(2.3.5)$$
By Proposition~(2.2), it is enough to show
$$V^0(\D_X[s]f^{s-\al})\subset M'\q\h{in}\,\,\, M.
\leqno(2.3.6)$$
Considering these mod $V^{>0}M$ and using the isomorphism (1.2.3), the assertion is further reduced to
$$\Gr_V^0(\D_X[s]f^{s-\al})\subset\Gr_V^0M'\q\h{in}\,\,\,\Gr_V^0M,
\leqno(2.3.7)$$
and follows from (2.3.4). This finishes the proof of Proposition~(2.3).
\bs\bs
\vbox{\centerline{\bf \S3. Proofs of the main theorems}
\bsn
In this section we prove the main theorems.}
\msn
{\bf 3.1.~Proof of Theorem~1.\,}
By Proposition~(2.1), we have to calculate
$$\D_X[\dd_t]f^{s-\al}\big/\D_X[\dd_t]f^{s-\al+1},$$
where we may assume $\al>0$ by the first inclusion in (1.1.5).
\sk
There is a natural isomorphism
$$\OO_X(*Z)[s]f^s\simto\OO_X(*Z)[s]f^{s-\al},
\leqno(3.1.1)$$
which is the identity on $\OO_X(*Z)$, and sends $s$ to $s-\al$.
(This is similar to the definition of the action of $t$ on $\D_X[s]f^s$ which sends $s$ to $s+1$, see \cite{Ma2}, \cite{Ka1}, etc.)
\sk
It induces isomorphisms of $\D_X$-submodules
$$\D_X[s]f^s\simto\D_X[s]f^{s-\al},
\leqno(3.1.2)$$
$$V^{\be}(\OO_X(*Z)[s]f^s)\simto V^{\be-\al}(\OO_X(*Z)[s]f^{s-\al})\q(\be\in\Q).
\leqno(3.1.3)$$
Note that these are not $\C[s]$-linear, and (3.1.2) does not hold if $\D_X[s]$ is replaced by $\D_X[\dd_t]$.
\sk
Combining these with Proposition~(2.3), we get
$$\aligned\Gr_V^0(\D_X[\dd_t]f^{s-\al})&=\Gr_V^{\al}(\D_X[s]f^s)+N\,\Gr_V^{\al}M_f\,\,(\subset\Gr_V^{\al}M_f)\\
&=G_j\Gr_V^{\be}M_f+N\,\Gr_V^{\be}M_f\,\,(\subset\Gr_V^{\be}M_f),\endaligned
\leqno(3.1.4)$$
where the last isomorphism is induced by the multiplication by $t^{-j}$ with $j:=\al-\be$, see the definition of $G$ in (1.3.1).
So Theorem~1 follows from the Verdier-type extension theorem for regular holonomic $\D$-modules as in (1.9). (Note that the composition of $(i_f)_*^{\!\D}$ and $\Gr_V^0$ is the identity for $\D_X$-modules supported in $Z$.) This finishes the proof of Theorem~1.
\msn
{\bf 3.2.~Proof of Theorem~2.\,} We may assume $\al\ne1$ by a remark after Question~1.
The first assertion of Theorem~2 follows from Theorem~1 and Proposition~1, since the two conditions
$$\Gr^G_j\M_f^{\be}\ne 0\q\h{and}\q\Gr^G_j\bl(\M_f^{\be}/N\M_f^{\be}\br)\ne 0$$
are equivalent to each other if $N$ is strictly compatible with the filtration $G$. Indeed, the strict compatibility implies the commutativity of $\Gr^G_j$ with ${\rm Coker}\,N$ (and $N$ is nilpotent). The last assertion of Theorem~2 then follows from (5). This finishes the proof of Theorem~2.
\msn
{\bf 3.3.~Proof of Theorem~3.\,} We may assume $\al\ne1$ as in the proof of Theorem~2 above.
By using (1.3.7) and Proposition~1, the assumption of Theorem~3 implies that there is
$$v\in G_j\M_f^{\be}\q\h{with}\q N^{\m_{\be}-1}v\notin G_{j-1}\M_f^{\be}.
\leqno(3.3.1)$$
We have to show
$$v\notin N\M_f^{\be}+G_{j-1}\M_f^{\be}.
\leqno(3.3.2)$$
If $v\in N\M_f^{\be}+G_{j-1}\M_f^{\be}$, then there is $v'\in\M_f^{\be}$ with
$$v=Nv'\mod G_{j-1}\M_f^{\be}.
\leqno(3.3.3)$$
Combined with (3.3.1), this implies
$$N^{\m_{\be}}v'\notin G_{j-1}\M_f^{\be}.
\leqno(3.3.4)$$
However, this contradicts the definition of $\m_{\be}$ in Theorem~3.
So (3.3.2) follows. This finishes the proof of Theorem~3.
\bs\bs
\vbox{\centerline{\bf \S4. Example}
\bsn
In this section we describe an example where Question~1 is negatively answered. We first show the following.}
\msn
{\bf 4.1.~Lemma.\,} {\it Let $h:=x^{14}+y^{14}-x^6y^6\in\C\{x,y\}$. Then the Milnor number $\mu_h$ is $141$, and the Jacobian ring $\C\{x,y\}/(\dd h)$ has a monomial basis consisting of $x^iy^j$ with $(i,j)$ contained in $I_0\cup I_1\cup{}^tI_1$, where ${}^tI_1$ is the transpose of $I_1$, and
$$I_0:=\bl\{(i,i)\mid i\in[0,10]\br\},\q
I_1:=\bl\{(i,j)\mid 1\les i-j\les 13,\,\,0\les j\les 4\br\}.$$
Moreover the spectral numbers $\{\al_{h,k}\}$ in $(1.7)$ are given by}
$$\bl\{\h{$\frac{j}{6}$}\mid j\in[1,11]\br\}\q\h{and}\q
\bl\{\h{$\frac{i}{14}+\frac{j}{6}$}\mid(i,j)\in[1,13]\times[1,5]\br\}\,\,\,(\h{with multiplicity $2$}).$$
\msn
{\it Proof.\,} We can calculate the Jacobian ring $\C\{x,y\}/(\dd h)$ by using the picture as below:
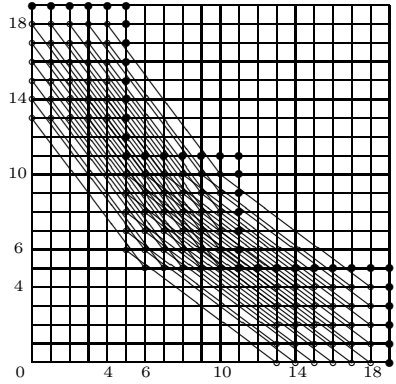
\sk
$$\setlength{\unitlength}{0.25cm}
\begin{picture}(20,20)
\multiput(1,1)(0,1){20}{\line(1,0){19}}
\multiput(1,1)(1,0){20}{\line(0,1){19}}
\multiput(1,14)(0,1){6}{\line(3,-4){6}}
\multiput(2,14)(0,1){6}{\line(3,-4){6}}
\multiput(3,14)(0,1){6}{\line(3,-4){6}}
\multiput(4,14)(0,1){6}{\line(3,-4){6}}
\multiput(5,14)(0,1){6}{\line(3,-4){6}}
\multiput(14,1)(1,0){6}{\line(-4,3){8}}
\multiput(14,2)(1,0){6}{\line(-4,3){8}}
\multiput(14,3)(1,0){6}{\line(-4,3){8}}
\multiput(14,4)(1,0){6}{\line(-4,3){8}}
\multiput(14,5)(1,0){6}{\line(-4,3){8}}
\multiput(1,14)(0,1){6}{\circle{.30}}
\multiput(2,14)(0,1){6}{\circle{.30}}
\multiput(3,14)(0,1){6}{\circle{.30}}
\multiput(4,14)(0,1){6}{\circle{.30}}
\multiput(5,14)(0,1){6}{\circle{.30}}
\multiput(14,1)(1,0){6}{\circle{.30}}
\multiput(14,2)(1,0){6}{\circle{.30}}
\multiput(14,3)(1,0){6}{\circle{.30}}
\multiput(14,4)(1,0){6}{\circle{.30}}
\multiput(14,5)(1,0){6}{\circle{.30}}
\multiput(11,6)(0,1){6}{\circle{.30}}
\multiput(10,6)(0,1){6}{\circle{.30}}
\multiput(9,6)(0,1){6}{\circle{.30}}
\multiput(8,6)(0,1){6}{\circle{.30}}
\multiput(7,6)(0,1){6}{\circle{.30}}
\multiput(6,7)(0,1){5}{\circle{.30}}
\multiput(12,6)(1,0){8}{\circle*{.45}}
\multiput(12,7)(0,1){6}{\circle*{.45}}
\multiput(7,12)(1,0){5}{\circle*{.45}}
\multiput(6,12)(0,1){8}{\circle*{.45}}
\multiput(20,1)(0,1){6}{\circle*{.45}}
\multiput(1,20)(1,0){6}{\circle*{.45}}
\put(0.4,0.2){\!$\scriptscriptstyle 0$}
\put(0.3,4.8){\!$\scriptscriptstyle 4$}
\put(0.3,6.8){\!$\scriptscriptstyle 6$}
\put(0,10.8){\!$\scriptscriptstyle 10$}
\put(0,14.8){\!$\scriptscriptstyle 14$}
\put(0,18.8){\!$\scriptscriptstyle 18$}
\put(4.8,0.2){$\scriptscriptstyle 4$}
\put(6.8,0.2){$\scriptscriptstyle 6$}
\put(10.6,0.2){$\scriptscriptstyle 10$}
\put(14.6,0.2){$\scriptscriptstyle 14$}
\put(18.6,0.2){$\scriptscriptstyle 18$}
\end{picture}$$
\skn
Here segments mean relations between the monomials corresponding to white vertices. Black vertices mean that the corresponding monomials are completely contained in the Jacobian ideal $(\dd h)\subset\C\{x,y\}$.
So the first assertion follows.
\sk
The Milnor number of $h$ is 141 by \cite{Ko}, and the above calculation is compatible with this.
Note that the above picture is closely related with a picture of the spectral numbers in the 2-dimensional non-degenerate Newton boundary case in \cite{Ar}. Here we have to shift by $(1,1)$ to adjust the difference between the Jacobian ring $\C\{x,y\}/(\dd h)$ and $\Om_h^2$ in (1.7).
\sk
We have an automorphism of $\Z^2$ defined by
$$\Z^2\ni(i,j)\mapsto(i+j,j)\in\Z^2.$$
This induces an bijection between a rectangle and a parallelogram with vertices
$$(0,0),(14,0),(0,6),(14,6)\q\h{and}\q(0,0),(14,0),(6,6),(20,6).$$
So the last assertion follows from \cite{St2} (since $n=2$).
This finishes the proof of Lemma~(4.1).
\msn
{\bf 4.2.~Example.\,} Let $\De\subset\C$ be a sufficiently small open disk. Set
$$h:=x^{14}+y^{14}-x^6y^6,\q g:=h+z^5,\q f_u:=g-u\,x^9y^2z\q(u\in\De).$$
These functions have non-degenerate Newton boundaries at 0 in the sense of Kouchinirenko (\cite{Ko}), and the $f_u\,\,\,(u\in\De)$ form a family of constant non-degenerate Newton boundary; in particular, they have the same Milnor number and the same monodromy, but they may have different Bernstein-Sato polynomials.
\sk
We have the Brieskorn lattices $H''_{f_u}$ and the Gauss-Manin systems $\G_{f_u}$ for $u\in\De$ as in (1.5).
They form the Brieskorn lattice $H''_{F,\De}$ and the Gauss-Manin system $\G_{F,\De}$ on $\De$ which are respectively locally free modules of rank $\mu$ over $\OO_{\De}\{\!\{\dd_t^{-1}\}\!\}$ and $\OO_{\De}\{\!\{\dd_t^{-1}\}\!\}[\dd_t]$, where
$$\OO_{\De}\{\!\{\dd_t^{-1}\}\!\}:=\bl\{\msum_{i\in\N}\,g_i\dd_t^{-i}\in\OO_{\De}[[\dd_t^{-1}]]\,\big|\,\msum_{i\in\N}\,g_it^i/i!\in\OO_{\De\times\C}|_{\De\times\{0\}}\br\},$$
and the $f_u$ ($u\in\De$) are identified with a function $F$ on $X\times\De$. \sk
Indeed, $H''_{F,\De,0}$ can be defined by
$$H''_{F,\De,0}:=\Om_{X\times\De/\De,0}^3/\ddd F\wedge\Om_{X\times\De/\De,0}^1,$$
(see also \cite{Gr}, \cite{SaK}, etc.)
The Gauss-Manin system $\G_{F,\De}$ can be defined by using $\delta(F-t)$ on $X\times\De\times\C$, see, for instance, \cite{bri}.
Note that the $\G_{f_u}$ are independent of $u$, and form a constant family, but the $H''_{f_u}$ do not.
The $V$-filtrations of Kashiwara and Malgrange on the $\G_{f_u}$ ($u\in\De$) define the $V$-filtration on $\G_{F,\De}$ with respect to the action of $t$ and $\dd_r$ as in (1.2) in a compatible way with the restriction to $\G_{F,\De}|_{\{u\}}=\G_{f_u}$ for $u\in\De$, since the $V$-filtrations on the $\G_{f_u}$ also form a constant family.
\sk
Set
$$\om_0:=\ddd x\wedge\ddd y\wedge\ddd z\in\Om_{X,0}^3,$$
where $(X,0):=(\C^3,0)$.
We will denote by $[\om]\in H''_{F,\De}$ the class of $\om\in\Om_{\C^3,0}^3$.
Define
$$v(\om):=\max\bl\{\al\in\Q\mid [\om]\in V^{\al}\G_{F,\De}\br\}.$$
This is closely related to the Newton filtration as in \cite{exp} (see also \cite{VK}).
For instance, we have
$$\h{$v(\om_0)=\frac{11}{30},\q v(x^9y^2z\,\om_0)=1+\frac{12}{30},\q v(x^{18}y^4z^2\om_0)=2+\frac{13}{30}$}.
\leqno(4.2.1)$$
Indeed, let $\Gamma_1$ be the face of the Newton polyhedron with vertices
$$(14,0,0),\,\,(6,6,0),\,\,(0,0,5).$$
It is contained in the hyperplane defined by
$$\ell_1(i,j,k):=\h{$\frac{1}{14}\,i+\frac{2}{21}\,j+\frac{1}{5}\,k=1$}.$$
Hence we get by \cite[Proposition 3.3]{exp} (which is applied to each $f_u$)
$$v(x^{i-1}y^{j-1}z^{k-1}\om_0)\ges\ell_1(i,j,k)=\h{$\frac{1}{14}\,i+\frac{2}{21}\,j+\frac{1}{5}\,k$},
\leqno(4.2.2)$$
that is,
$$[x^{i-1}y^{j-1}z^{k-1}\om_0]_{f_u}\in V^{\al}\G_{f_u}\,\,\,\,\,(u\in\De),$$
with $\al=\ell_1(i,j,k)$ if $(i,j,k)\in\C^3$ is contained in the cone of the face $\Gamma_1$.
\sk
For the proof of (4.2.1), we may then restrict to $u=0$, and the assertion is reduced to the calculation of the Jacobian ring $\C\{x,y,z\}/(\dd g)$ (since $f_0=g$). This is further reduced to that for $\C\{x,y\}/(\dd h)$ by using the Thom-Sebastiani type theorem for Brieskorn lattices (see \cite{Ma1}, \cite{SS}). Then (4.2.1) follows from Lemma~(4.1).
\sk
On the other hand, we have
$$\aligned\dd_u[\om_0]&=\dd_t[x^9y^2z\,\om_0]\in V^{12/30}\G_{F,\De},\\
\dd_u\dd_t[x^9y^2z\,\om_0]&=\dd_t^2[x^{18}y^4z^2\om_0]\in V^{13/30}\G_{F,\De},\\
\dd_u\dd_t^2[x^{18}y^4z^2\om_0]&=\dd_t^3[x^{27}y^6z^3\om_0]\in V^{14/30}\G_{F,\De}.\endaligned
\leqno(4.2.3)$$
This can be verified by using the relation
$$\dd_u\delta(F-t)=-(\dd_uF)\dd_t\delta(F-t),$$
(see for instance \cite{bri}, also \cite{Gr}, \cite{SaK}, etc.) By a similar argument we get
$$\h{$\Gr_V^{\al}[\om]\in\Gr_V^{\al}\G_{F,\De}$ is constant on $\De$ for $\om\in\Om_{X,0}^3$ and $\al=v(\om)$.}
\leqno(4.2.4)$$
\sk
These calculations imply that $-13/30$ is a {\it simple} root of $b_{f_u}(s)$ for $u\in\De^*$ by using Malgrange's formula in (1.5) together with the Taylor expansion of $[\om_0]\in\G_{F,\De,0}$ in $u$.
Indeed, for the constant part $\G_{F,\De}|_{u=0}=\G_g$, we have the decomposition
$$\aligned\G_g&=\widehat{\mopl}_{\al\in\Q}\,\G_f^{\al},\\
\h{with}\q\G_g^{\al}&={\rm Ker}(\dd_tt-\al)^2\subset\G_f,\endaligned$$
(see for instance \cite{bri}). So we get the expansion
$$[\om_0]=\msum_{r\in\N}\,u^r\xi_r=\msum_{r\in\N,\al\in\Q}\,u^r\xi_r^{\al}\q\h{in}\,\,\,\G_{F,\De},$$
with $\xi_r\in\G_g$, $\xi_r^{\al}\in\G_g^{\al}$, and we have for $r=0,1,2,\dots$
$$\xi_r=\msum_{\al\in\Q}\,\xi_r^{\al}=\dd_t^r[(x^9y^2z)^r\om_r]_g\q\h{in}\,\,\,\G_g.$$
\sk
Set
$$\be_0=\h{$\frac{13}{30}$}.$$
We have
$$\xi_r^{\be_0}=0\q\h{unless}\,\,\,\,r=2.
\leqno(4.2.5)$$
Indeed, $\xi_r\in V^{>\be_0}\G_g$ for $r>2$, and
$$\xi_r^{\al}=0\q\h{if}\,\,\,\,\al\notin\h{$\frac{11+r}{30}+\frac{\N}{14}+\frac{\N}{6}$}\q(r=0,1),$$
by using (4.2.1), Lemma~(4.1) together with the Thom-Sebastiani type theorem for Gauss-Manin systems
$$\G_h\otimes_K\G_{z^5}\simto\G_g,$$
where $K:=\C\{\!\{\dd_t^{-1}\}\!\}[\dd_t]$, see (1.6.2).
Note that $[z^r\ddd z]_{z^5}$ is annihilated by $\dd_tt-\al$ so that
$$[z^r\ddd z]_{z^5}\in\G_{z^5}^{\al}\q\h{with}\,\,\,\,\al=(r+1)/5.$$
(However, it is rather nontrivial to determine the expansion of $[\om]_h\in\G_h$ even if $\om\in\Om_{\C^2,0}^2$ is a monomial in the case where $h$ is not a weighted homogeneous polynomial.)
\sk
We have moreover
$$\xi_2^{\be_0}=\Gr_V^{\be_0}\dd_t^2[x^{18}y^4z^2\om_0]_g\in N\Gr_V^{\be_0}H''_g\subset\Gr_V^{\be_0}\G_g.
\leqno(4.2.6)$$
(However this does not hold for $\Gr_V^{\be_0}\dd_t^2[x^{10}y^{10}z^2\om_0]_g$ as one might expect from the picture of spectral numbers in \cite{Ar}.)
Indeed, (4.2.6) follows from
$$\bl(\dd_tt-\h{$\frac{43}{30}$}\br)[x^4y^4z^2\om_0]_g=c_1\,\dd_t[x^{18}y^4z^2\om_0]_g\q\q(c_1\in\C^*),$$
or equivalently
$$\bl(\dd_t^{-1}-\h{$\frac{30}{43}$}t\br)[x^4y^4z^2\om_0]_g=c'_1\,[x^{18}y^4z^2\om_0]_g\q\q(c'_1\in\C^*).$$
Here $\eta$ in (1.5.1) can be given by the contraction of $\om=x^4y^4z^2\om_0$ with the vector field
$$c_2^{-1}\bl(\h{$\frac{2}{21}x\dd_x+\frac{1}{14}y\dd_y+\frac{1}{5}z\dd_z$}\br)\q\q\bl(c_2=v(x^4y^4z^2\om_0)=\h{$\frac{43}{30}$}\br),$$
and the above equality follows.
\sk
To simplify the notation, set in the notation of (6)
$$\G^{\be_0}:=\Gr_V^{\be_0}\G_{f_u}=H^2(F_{\!f_u,0},\C)_{\ee(-\be_0)}.$$
Combining (4.2.3) and (4.2.5-6), we get for $u\in\De^*$
$$0\ne G_0\,\G^{\be_0}=\Gr_V^{\be_0}\Ht''_{f_u}\subset N\G^{\be_0}\subset\G^{\be_0},
\leqno(4.2.7)$$
where the filtration $G$ on $\G^{\be_0}$ is defined by the same way as in (1.3).
\sk
Indeed, the second spectral number $\al_{g,2}$ satisfies the relation
$$\al_{g,2}=\h{$\frac{11}{30}+\frac{1}{14}>\frac{13}{30}$}=\be_0,
\leqno(4.2.8)$$
by using Lemma~(4.1) and (1.7.2).
So it is enough to calculate the Taylor expansion of $[\om_0]\in H''_{F,\De}$ in order to determine $G_0\,\G^{\be_0}=\Gr_V^{\be_0}\Ht''_{f_u}$. So (4.2.7) follows.
\sk
Since
$$\dim\G^{\be_0}=2\q\h{with}\q N^2=0,$$
(by using \cite{ST}), we thus get
$$\C\cong G_0\,\G^{\be_0}=N\G^{\be_0}\ne\G^{\be_0}\cong\C^2.
\leqno(4.2.9)$$
In the notation of the introduction, this implies
$$\Gr^G_0M_{f_u}^{\be_0}\ne 0,\q\Gr^G_0(M_{f_u}^{\be_0}/NM_{f_u}^{\be_0})=0.
\leqno(4.2.10)$$
\msn
{\bf 4.3.~Remark.\,} The argument in Example~(4.2) can be generalized to the case where
$$f_u=x^{2a}+y^{2a}-x^{2b}y^{2b}+z^c-u\,x^{a+b-1}y^{b-1}z,$$
with $a,b,c$ mutually prime positive integers satisfying
$$a>2b>c,\q\h{$\frac{1}{2a}>\frac{2}{c}-\frac{1}{b}$}\,.
\leqno(4.3.1)$$
Indeed, the hyperplane passing through
$$(2a,0,0),\,\,(2b,2b,0),\,\,(0,0,c)$$
is defined by
$$\ell_1(i,j,k):=\h{$\frac{1}{2a}\,i+\frac{a-b}{2ab}\,j+\frac{1}{c}\,k=1$},$$
and we have
$$\aligned\ell_1(1,1,1)&=\h{$\frac{1}{2b}+\frac{1}{c}$},\\
\ell_1(a+b,b,2)&=\h{$1+\frac{2}{c}$},\\
\ell_1(2a+2b-1,2b-1,3)&=\h{$2-\frac{1}{2b}+\frac{3}{c}$}.\endaligned$$
\sk
Note that $x^{2a+2b-2}y^{2b-2}$ is identified with $x^{4b-2}y^{4b-2}$ up to a nonzero constant multiple in the Jacobian ring of $h:=x^{2a}+y^{2a}-x^{2b}y^{2b}$.
The last condition of (4.3.1) is needed to prove the generalizations of (4.2.5) and (4.2.8) with $\be_0:=\frac{3}{c}-\frac{1}{2b}$.
If we set $c=2b-1$, then (4.3.1) becomes the condition
$$b(2b-1)>2a>4b,$$
in particular $b\ges 3$.
\msn
{\bf 4.4.~Remark.\,} The calculation in Example~(4.2) implies that the Milnor monodromy $T$ on the vanishing cohomology $H:=H^{n-1}(F_{\!f,0},\C)$ and $e^{-2\pi i\dd_tt}$ on $\Ht''_f/t\Ht''_f$ do not have the same Jordan normal form, where $f=f_u$ for $u\in\De^*$. (This solves an old open problem in \cite{sup}.)
\sk
Indeed, the induced filtration $V$ on the $\mu$-dimensional $\C$-vector space $\Ht''_f/t\Ht''_f$ splits by the eigenvalues of the action of $s:=-\dd_tt$, and we have the canonical isomorphisms
$$\aligned\mopl_{\al}\,\bl(\Gr_V^{\al}(\Ht''_f/t\Ht''_f),\,e^{-2\pi i\dd_tt}\br)&=\mopl_{\be,j}\,\bl(\Gr^G_j\Gr_V^{\be}\G_f,\,e^{-2\pi i\dd_tt}\br)\\
&=\mopl_{\la,j}\,\bl(\Gr^G_jH_{\la},\,T\br),\endaligned$$
where $(H_{\la},G)$ is identified with $(\Gr_V^{\be}\G_f,G)$ for $\la=\exp(-2\pi i\be)$.
The problem is thus reduced to the comparison between the actions of $N$ on $H_{\la}$ and $\mopl_j\,\Gr^G_jH_{\la}$. So the assertion follows from (4.2.9).

\end{document}